\documentclass[12pt]{amsart}

\usepackage{fancyhdr}

\newcommand{\pic}{{\rm Pic}}
\newcommand{\bs}{{\rm Bs}}
\newcommand{\bl}{\bs\mid L\mid} 
\newcommand{\be}{\bs\mid 2L\mid} 
\newcommand{\p}{\mathbb{P}}
\newcommand{\co}{\mathbb{C}}
\newcommand{\z}{\mathbb{Z}}
\newcommand{\rt}{\rightarrow}

\theoremstyle{plain}
\newtheorem{lemma}{Lemma}[section]
\newtheorem{teo}[lemma]{Theorem}

\newtheorem{propo}[lemma]{Proposition}
\newtheorem{coro}[lemma]{Corollary}
\newtheorem{con}[lemma]{Conjecture}

\theoremstyle{definition}
\newtheorem{defi}[lemma]{Definition}

\newtheorem{example}[lemma]{Example}

\theoremstyle{remark}

\begin{document}
\setlength{\baselineskip}{15pt}

\author{Antonio Laface}
\title{A very ampleness result}
\email{laface@mat.unimi.it}
\maketitle

\begin{abstract}
Let $(M,L)$ be a polarized manifold. The aim of this paper is to 
establish a connection between the generators of the graded algebra
$\bigoplus_{i\geq 1}H^0(M,iL)$ and the very ampleness of the line bundle $rL$.
Some applications are given.
\end{abstract}

\section{Introduction}
Let $L$ be an ample line bundle on an algebraic manifold $M$, the problem 
of finding the least $n$ such that $nL$ is very ample is a basic one in 
the classification theory of polarized varieties. Many attempts were made 
in order to establish a general formula. A fundamental result due to 
Matsusaka \cite{ma} says that there is a constant $c$, depending only 
on the Hilbert polynomial of $(M,L)$, such that $cL$ is very ample. 
Moving in another direction Fujita conjectured \cite[\S2, Conjecture b]{f4}
that $K_M+(n+2)L$ is very
ample and $K_M+(n+1)L$ is spanned for every polarized manifold $(M,L)$ with
$n=\dim M$.  In this paper we consider the graded algebra:
\[
G(M,L)=\bigoplus_{i\geq 1}H^0(M,iL).
\]
Now the generators of this algebra describe the whole structure of $(M,L)$. 
They also allows us to determine the embedding of 
$M$ via the linear system $\mid rL\mid$ for $r\gg 0$, so there must be a 
connection between these generators and the very ampleness of $rL$.
The aim of this paper is pointing out this conection.
The paper is organized as follows: in section 2 we prove the main theorem 
asserting the very ampleness of $rL$, outside the base locus of $L$, for $L$ 
a $r$-generated line bundle. This also implies the very ampleness of $2L$
for
any ample $2$-generated line bundle \ref{mt2}.
In section 3 some applications of the
main theorem \ref{tm} are given, concerning:
a) sectionally hyperelliptic polarized varieties of type $(-)$, b) surfaces of
general type, c) del Pezzo surfaces.
\section{The main theorem} All notation used in this paper are standard in
algebraic geometry.
Let $L$ be a line bundle on a projective manifold $M$ and consider the 
associated graded algebra $G(M,L)$.
The following definition generalizes the notion of simply generated line 
bundle.
\begin{defi}
Let $r\geq 1$ be an integer. A line bundle $L$ on $M$ is $r$-gen\-er\-ated if
$G(M,L)$ is generated by the sections of $H^0(M,L),\ldots ,H^0(M,rL)$.
We also say that the pair $(M,L)$ is $r$-generated.
\end{defi}
There is a connection between the $r$-generation of $L$ and the very 
ampleness of the line bundle $rL$; this connection, extending a known fact holding for simply generated line bundles, is expressed by the following:
\begin{teo}
Given a polarized manifold $(M,L)$ with $L$ effective and
$r$-generated, then $\varphi_{\mid rL\mid}$ is an embedding of
$M\backslash\bl$.
\label{tm}
\end{teo}
\begin{proof}
Set $N=M\backslash\bl$ and $\varphi=\varphi_{\mid rL\mid}$.
First of all we observe that $\varphi$ is well defined in $N$, since
$L$ is spanned on $N$ and hence so are its multiples.\\
By hypothesis there is an integer $k\geq 1$ such that $kL$ is very ample.
Take $p,q\in N$ and by contradiction suppose that no section of $H^0(iL)$ with
$1\leq i\leq r$ separates  $p$ and $q$. Then each section $s\in H^0(M,iL)$ 
that vanishing at $p$ must vanish at $q$ too.\\
Consider two sections $s_1,s_2\in H^0(iL)$, and define $\gamma,\lambda$ so 
that $s_1(p)=\gamma s_1(q)$ e $s_2(p)=\lambda s_2(q)$.
If we take $s=s_1(p)s_2-s_2(p)s_1$ we have that $s(p)=0$ and necessarily
$0=s(q)=s_1(q)s_2(q)(\gamma-\lambda)$. This implies that $\gamma=\lambda$. 
we note that if one of the two sections vanishes at $q$ it must vanish also 
at $p$, otherwise it will separate $p$ and $q$, but in this situation we may 
equally take $\gamma=\lambda$.\\
The previous argument shows that there exist some constants $\gamma_i$ 
associate to each 
$H^0(iL)$ for $1\leq i\leq r$ such that for each section 
$s\in H^0(iL)$ we have $s(p)=\gamma_i s(q)$.\\
Now note that since $G(M,L)$ is an algebra,
the constants $\gamma_i$ have to satisfy the relations:
\[
\gamma_i=\gamma^i \  {\rm where}\  \gamma:=\gamma_1.
\]
For, since $L$ is effective there exists $\sigma\in H^0(L)$, then
$\sigma^i(p)=\gamma^i\sigma^i(q)$, but $\sigma^i\in H^0(iL)$ so we have also
$\sigma^i(p)=\gamma_i\sigma^i(q)$ and this implies the previous equation.\\
Now let us consider a section $s\in H^0(kL)$.
By the $r$-generation hypothesis the sections of $H^0(M,kL)$ are
linear combinations of products of sections in $H^0(iL)$ with
$1\leq i\leq r$. 
Then we can write $s=\sum \xi_j s_j$ where $\xi_j\in\co$ and
$s_j=\prod_{i=1}^r\prod_{n_i=1}^{\alpha_i}\sigma_{n_i}$
with $\sigma_{n_i}\in H^0(iL)$, where $\alpha_i$ is the number of sections
of $H^0(iL)$ that appears in the product. Note that $\sum_{i=0}^r i\alpha_i=k$ because $s_j\in H^0(kL)$.
We have
\[
s(p)=\sum\xi_j s_j(p)=\sum\xi_j 
(\prod_{i=1}^r\prod_{n_i=1}^{\alpha_i}\sigma_{n_i}(p))=
\sum\xi_j (\prod_{i=1}^r\prod_{n_i=1}^{\alpha_i}\gamma^i\sigma_{n_i}(q))=
\]
\[
=\sum\xi_j (\prod_{i=1}^r(\gamma^{i\alpha_i}\prod_{n_i=1}^{\alpha_i}
\sigma_{n_i}(q)))=
\sum\xi_j\gamma^{({\sum_{i=1}^r i\alpha_i})}(\prod_{i=1}^r\prod_{n_i=1}^
{\alpha_i}\sigma_{n_i}(q))=
\]
\[
=\sum\xi_j\gamma^k (\prod_{i=1}^r\prod_{n_i=1}^{\alpha_i}\sigma_{n_i}(q))=
\gamma^k\sum\xi_j s_j(q)=\gamma^k s(q).
\]
We thus obtained that the sections of $H^0(kL)$ do not separate points, 
which is absurd.
Then there exists a section $\sigma\in H^0(iL)$ wich separates
$p$ and $q$, i.e. $\sigma(p)=0$ and $\sigma(q)\neq 0$.
Since $L$ is spanned outside $\bl$ we may take a section $\delta\in H^0(L)$ such that $\delta(q)\neq 0$. 
Now the section $\sigma\delta^{r-i} H^0(rL)$ separates $p$ and $q$.\\
In order to show that the map $\varphi$ is an immersion at each point
$p\in N$ we have to show that for each vector $\tau\in T_pM$,
there exists a section $\sigma\in H^0(iL)$ such that $d\sigma(\tau)\neq 0$
and $\sigma(p)=0$.
By contraddiction suppose that such a section does not exist, then for each 
$1\leq i\leq r$ there exists a costant $\eta_i$ such that for each
 $s\in H^0(iL)$, $ds(\tau)=\eta_i s(p)$.
To show this let $s_1,s_2\in H^0(iL)$ and consider $\alpha,\beta$ such that $ds_1(\tau)=\alpha s_1(p)$
and $ds_2(\tau)=\beta s_2(p)$. Now consider the section $s=s_1(p)s_2-s_2(p)s_1$;
we have $s(p)=0$ and $0=ds(\tau)=s_1(p)s_2(p)(\beta-\alpha)$. This implies 
$\alpha=\beta$.\\
Define $\eta:=\eta_1$. Take $\sigma^i\in H^0(iL)$ then $\eta_i \sigma^i(p)=d[\sigma^i](\tau)=i d\sigma({\tau})\sigma^{i-1}(p)=i \eta\sigma^i(p)$.
Hence we have:
\[
\eta_i=i\eta.
\]
Now for a section $s\in H^0(kL)$ we have:
\[
ds(\tau)=\sum\xi_j ds_j(\tau)=\sum\xi_j 
d[\prod_{i=1}^r\prod_{n_i=1}^{\alpha_i}\sigma_{n_i}](\tau)=
\]
\[
=\sum\xi_j (\sum_{h=1}^r d[\prod_{n_h=1}^{\alpha_h}
\sigma_{n_h}](\tau)\prod_{i=1,i\neq h}^r\prod_{n_i=1}^{\alpha_i}
\sigma_{n_i}(p))=
\]
\[
=\sum\xi_j (\sum_{h=1}^r(\sum_{m=0}^{\alpha_h}d\sigma_m(\tau)
\prod_{n_h=1,n_h\neq m}^{\alpha_h}\sigma_{n_h}(p))
\prod_{i=1,i\neq h}^r\prod_{n_i=1}^{\alpha_i}\sigma_{n_i}(p))=
\]
\[
=\sum\xi_j (\sum_{h=1}^r(\sum_{m=0}^{\alpha_h}h\eta_1
\prod_{n_h=1}^{\alpha_h}\sigma_{n_h}(p))
\prod_{i=1,i\neq h}^r\prod_{n_i=1}^{\alpha_i}\sigma_{n_i}(p))=
\]
\[
=\sum\xi_j (\sum_{h=1}^r h\alpha_h\eta_1)
\prod_{i=1}^r\prod_{n_i=1}^{\alpha_i}\sigma_{n_i}(p)=
\sum\xi_j k\eta_1
\prod_{i=1}^r\prod_{n_i=1}^{\alpha_i}\sigma_{n_i}(p)=\eta_k s(p).
\]
This implies that the sections of $H^0(kL)$ do not separate  $p$ from the vector
$\tau$, but this contradicts the very ampleness of $kL$.
\end{proof}

Theorem \ref{tm} immediately gives the following
\begin{coro}
Let $(M,L)$ a polarized manifold with $L$ spanned and $r$-ge\-ne\-ra\-ted, then $rL$ is very ample.
\end{coro}
Note that theorem \ref{tm} cannot be inverted. The following examples 
show very ample line bundles $L$ which are not $1$-generated.
\begin{example}
{\rm Let $(S,L)$ be an abelian surface polarized by a very ample line bundle. 
Let $L^2=2d$, since $S$ is abelian $K_S$ is the trivial bundle, so we have 
$h^i(L)=h^i(K_S+L)=0$ for $i=1,2$. It follows that 
$h^0(L)=\chi(L)=\chi ({\mathcal O}_S)+\frac{L^2-LK_S}{2}=d$.\\
Let $S^2(H^0(L))$ denote the second symmetric power of $H^0(L)$. Then:
\[
\dim S^2(H^0(L))=\frac{d(d+1)}{2}.
\]
Now if $L$ is simply generated we have $S^2(H^0(L))\cong H^0(2L)$, but
$h^0(2L)=\chi ({\mathcal O}_S)+\frac{4L^2-2LK_S}{2}=4d$. So if $d\leq 6$
we have $h^0(2L)=4d>\frac{d(d+1)}{2}=S^2(H^0(L))$. 
It is well known that there exist abelian surfaces of degree $2d=10$ in
$\p^4$.}
\end{example}
\begin{example}
{\rm 
Let $Q=\p^1\times\p^1\subset\p^3$, we have $\pic (Q)\cong\z\times\z$.
Let $L_1,L_2\in\pic (Q)$ be two generators; then $2L_1+4L_2$ is very 
ample [\cite{ha}, ch. 2.18, p. 380]. So there exists a smooth curve 
$C\in\mid 2L_1+4L_2\mid$.
Now consider the polarized curve $(C,L)$, where 
$L={\mathcal O}_{\p^3}(1)_{C}$. Since $C$ is smooth we have:
$2g(C)-2=(2L_1+4L_2+K_Q)(2L_1+4L_2)=(2L_1+4L_2-2L_1-2L_2)(2L_1+4L_2)=
2L_2(2L_1+4L_2)=4$. Hence $g(C)=3$ and $\deg L=(L_1+L_2)(2L_1+4L_2)=6>2g(C)-1$.
Then $h^1(L)=0$ and $h^0(L)=4$.
From this it follows that:
$
\dim S^2(H^0(L))=h^0(2L)=10.
$
Now consider the exact sequence:
\[
0\rt {\mathcal I}_C\rt {\mathcal O}_{\p^3}\rt {\mathcal O}_C\rt 0
\]
Tensoring with ${\mathcal O}_{\p^3}(2)$ and taking cohomology, we obtain:
\[
0\rt H^0({\mathcal I}_C(2))\rt H^0({\mathcal O}_{\p^3}(2))\rt H^0(2L)\rt\ldots.
\]
Since $C$ is contained in a quadric it follows that 
$H^0({\mathcal I}_C(2))\neq 0$,
and recalling that $H^0({\mathcal O}_{\p^3}(2))=S^2(H^0(L))$ we immediately
see that the map: $S^2(H^0(L))\rt H^0(2L)$ cannot be surjective.
}
\end{example}
Theorem \ref{tm} can be regarded as a generalization of the following well 
known fact:
\begin{propo}
Let $(M,L)$ be a polarized manifold and assume that $L$ is simply generated 
(i.e. $1$-generated), then $L$ is very ample.
\end{propo}
We have only to show that in this case $\bl=\emptyset$. By contradiction 
suppose that 
there exists $p\in\bl$ and let $tL$ be very ample. since $L$ is $1$-generated, 
each section of $tL$ is a sum of products of sections belonging to $H^0(L)$, 
then $tL$ would not be spanned at $p$.\\
Theorem \ref{tm} also allows us to prove the following:
\begin{teo}
\label{mt2}
Let $(M,L)$ be a polarized manifold with effective and $2$-ge\-ne\-ra\-ted line
 bundle; then $2L$ is very ample.
\end{teo}
\begin{proof}
From teorem \ref{tm} we have only to show that:
\begin{itemize}
\item $2L$ is spanned.
\item $2L$ separates $p,q$, with at least one in $\bl$.
\item $2L$ separates $p\in\bl$ from each $\tau\in T_pM$.
\end{itemize}
Let $tL, t\geq 2$, be a very ample line bundle. Each section $s\in H^0(tL)$ 
is of the form
\begin{equation}
\label{fo}
s=\sum(\prod_{i=0}^a\alpha_i\prod_{j=0}^b\beta_j),
\end{equation}
with $\alpha_i\in H^0(2L)$ and $\beta_j\in H^0(L)$ so that $2a+b=t$.
Now $2L$ is spanned outside $\bl$, since $\be\subset\bl$.
Let $p\in\bl$, if every $\alpha\in H^0(2L)$ vanishes at $p$ then each product in \ref{fo}
must vanish, so $tL$ is not spanned.\\
If $p\in\bl$ and $q\not\in\bl$ take $\beta\in H^0(L)$ such that $\beta(q)\neq 
0$, $\beta^2\in H^0(2L)$ separates $p$ and $q$.
Let $p,q\in\bl$. As in the proof of teorem \ref{tm} the same shows that 
there exists a section in $H^0(L)$ or in $H^0(2L)$ separating $p$ from $q$, 
then $2L$ separates $p$ and $q$.\\
If $p\in\bl$ and $\tau\in T_pM$, then as before we know that exists a section 
of $H^0(L)$ or of $H^0(2L)$ separating $p$ and $\tau$. But in this case we may have a $\beta\in H^0(L)$ with
$d\beta_p(\tau)\neq 0$. This doesn't tell anything on $H^0(2L)$ because $\beta^2(p)=0$ but also
$d\beta^2_p(\tau)=0$. Now we prove that, if $2L$ does not separate $p$ and $\tau$, such a $\beta$ does not exists. We have two cases:\\
$t$ is even. Let $\eta$ be as in theorem \ref{tm} such that 
$d\alpha_p(\tau)=\eta\alpha(p)$, for each $\alpha\in H^0(2L)$.
Let $t=2r$ and consider $s\in H^0=(tL)$:
\[
ds_p(\tau)=\sum(d(\prod_{i=1}^a\alpha_i\prod_{2j=0}^b\beta_j)_p(\tau))=
\sum(d(\prod_{i=1}^r\alpha_i)_p(\tau))=
\]
\[
\sum\sum_{i=1}^r(\prod_{j\neq i}\alpha_j(p))d{\alpha_i}_p(\tau)=
\sum\sum_{i=1}^r\eta(\prod_{j=1}^r\alpha_j(p))=
\]
\[
\sum r\eta(\prod_{j=1}^r\alpha_j(p))=r\eta s(p)
\] 
So $tL$ does not separate $p$ and $\tau$.\\
$t$ is odd. Since $\beta_j(p)=0$ for each $j$, we have:
\[
s(p)=\sum(\prod_{i=0}^a\alpha_i(p)\prod_{2j+1=0}^b\beta_j(p))=0.
\]
Then $tL$ would not span in $p$.
\end{proof}
\section{Some applications}
\noindent
In this section we discuss some applications of theorem \ref{tm}.\\
First we consider a polarized manifold $(M,L)$ with 
$\bs\mid L\mid=\{p_1,\ldots, p_n\}$ (a finite set) and
we obtain some information on the generator of $G(M,L)$.
To apply the main technique: "the Apollonius method" see \cite{f}, recall that
given $(M,L)$ be a polarized manifold,
a sequence of irreducible and reduced subvarieties $D_i$ of $M$ such that
$D_i\in\mid L_{D_{i-1}}\mid$ and $D_1\in\mid L\mid$ is called a 
{\bf ladder} of $(M,L)$ and the $D_i$ are said the {\bf rungs} of the
ladder.
In the second application we study some connection between the generator of the
canonical ring $G(S,K_S)$ of a surface of general type and the very 
ampleness of $tK_S$. 
Finally we study the very ampleness of $-tK_S$ for a del Pezzo surface.

\subsection{Hyperelliptic polarized manifold of type $(-)$}

In this subsection we apply the main theorem 
to provide an evidence for the first conjecture of Fujita in a special case.
First of all we recall Fujita's conjecture, see 
\cite[\S2, Conjecture b]{f4} 
\begin{con}
For every polarized manifold $(M,L)$ with $\dim M=n$, $K_M+(n+2)L$ is very ample.
\end{con}
We recall that the delta-genus $\Delta(M,L)$ of a polarized manifold $(M,L)$ is
defined as $\Delta(M,L)=L^n+n-h^0(M,L).$
Here we deal with polarized manifolds with $\Delta(M,L)=1$. In 
this case for $L^n\geq 3$ it is known \cite{f} that $L$ is very ample. It is also easy 
to prove that for $L^n=2$ the line bundle $(g+1)L$ is very ample \cite{lps1}.\\
The case $L^n=1$ is very difficult to classify and there are only partial
results, see \cite{f3}. Here we consider only polarized manifolds sectionally
hyperelliptic of type $(-)$.
In this case Fujita \cite{f3} proved the following facts:
\begin{equation} 
H^q(M,tL)=0 \ {\rm for} \ q\geq 1, \ t\in\z, \label{pi}
\end{equation}
\begin{equation} 
K_M=(2g-n-1)L \ {\rm where} \ g=g(M,L) \ {\rm and} \ n=\dim M,
\label{ad}
\end{equation}
\begin{equation} \pic M\cong\z, \ {\rm generated \ by} \ L \end{equation}
\begin{equation} (M,L) \ {\rm has \ a \ ladder}.\label{lp} \end{equation}
Properties \ref{pi} and \ref{ad} are true also for the
rungs $(D,L_D)$ of this ladder.
In this situation we may apply theorem \ref{tm} to prove the very ampleness of $K_M+(n+2)L=(2g+1)L$.\\
We start with a lemma on the $r$-generation inspired by \cite[2.3]{f} 
\begin{lemma}
Let $(M,L)$ be a polarized manifold with $\Delta(M,L)=d(M,L)=1$, sectionally 
hyperelliptic of type $(-)$; if the line bundle $L_C$ is $r$-generated, then $L$ is $r$-generated.
\end{lemma}
\begin{proof}
Consider a ladder of $(M,L)$: $(M_1,L_1)=(M,L),\ldots,(M_n,L_n)=(C,L_C)$ and proceed by induction on $n$.\\
For $n=1$ the assumption is true. Now suppose that it is true for $k+1$ and 
consider the exact sequence:
\[
0\rt {\mathcal O}_{M_k}(-L_k)\rt {\mathcal O}_{M_k}\rt {\mathcal O}_{M_{k+1}}\rt 0,
\]
tensor with ${\mathcal O}_{M_k}(tL_k)$ and consider the exact cohomology 
sequence, recalling \ref{pi} we get:
\begin{equation}
\label{f1}
0\rt H^0(M_k,(t-1)L_k)\rt H^0(M_k,tL_k)\rt H^0(M_{k+1},L_{k+1})\rt 0. 
\end{equation}
So for each $t$ the restriction map $\psi$:
\[
H^0(M_k,tL_k)\rt H^0(M_{k+1},tL_{k+1})
\]
is surjective. Let us take a set of generators $\gamma_1,\ldots,\gamma_i$ of 
$G(M_{k+1},L_{k+1})$ and their 
inverse images $\eta_1,\ldots,\eta_i$ in $G(M_k,L_k)$, if $\sigma\in H^0(M_k,L_k)$ is the defining section of $M_{k+1}$ then $G(M_k,L_k)$ is generated by the $\eta_j$ and $\sigma$.
To show this, observe that by \ref{f1} $H^0(M_k,L_k)$ is generated by 
$\sigma$ and by the  $\eta_j$ that belong to $H^0(M_{k+1},L_{k+1})$. 
Proceeding by induction on $t$ it is simple to prove the assertion.
Finally observe that $\eta_1,\ldots,\eta_i$ and $\sigma$ belong to 
$H^0(M_k,tL_k)$ with $t\leq r$. This completes the proof.
\end{proof}

Considering the ladder of $(M,L)$, by adjunction from \ref{ad} we have:
\[
K_{M_i}=(2g-i-1)L_i.
\]
Then $K_C=(2g-2)L_C=(2g-2)p$, with $p=\bl$. Hence $g(M,L)=g=g(C)$ the genus of the curve $C$. Then since $C$ is an hyperelliptic curve we have that $L_C$ is exactly $2g+1$-generated. Then by the lemma $(M,L)$ is $2g+1$-generated.\\ 
The following theorem extends what proved in \cite{lps2} for $g=1$.
\begin{teo}
Let $(M,L)$ be a polarized manifold with $\Delta(M,L)=L^n=1$, sectionally 
hyperelliptic of type $(-)$, then $(2g+1)L$ is very ample.
\label{t5}
\end{teo}
\begin{proof}
Since $L$ is $(2g+1)$-generated, by \ref{tm} we know that $\varphi_{\mid (2g+1)L\mid}$ 
is an embedding of $M\backslash\{p\}$.\\
To see that $(2g+1)L$ is spanned on $p$ we proceed by induction.
Clearly $(2g+1)L_C$ is very ample.
Consider a smooth $D\in\mid L\mid$ (we know that such a $D$ exists by the ladder property),
then consider the exact sequence:
\[
0\rt H^0(2gL)\rt H^0((2g+1)L)\rt H^0((2g+1)L_D)\rt 0.
\]
By induction we have that  $(2g+1)L_D$ is spanned on $p$ then there exists $\sigma\in H^0((2g+1)L_D)$ with $\sigma(p)\neq 0$. By the suriectivity of 
the restriction map we have a $\gamma\in H^0((2g+1)L)$ such that $\gamma(p)\neq 0$.\\
To see that $(2g+1)L$ defines an embedding in $p$, consider a vector 
$\tau\in T_pM$. We say that exists a section $\eta\in H^0(M,L)$ such that $d\eta (\tau)\neq 0$. To construct $\eta$ we may take an hypersurface $D\in\mid L\mid$ such that $\tau\not\in T_pD$ then the section that define $D$ clearly
have a differential that does not vanish on $\tau$. Suppose by contraddiction that for each 
$D\in\mid L\mid$ $\tau\in T_pD$ then the intersection of $n$ such $D_i$ cannot be transverse, so $L^n>1$, absurd. Section $\eta^{2g+1}$ defines an embedding in $p$.\\
Now consider a point $q\in M$ different from $p$, we have to find a section in $H^0((2g+1)L)$ that separates $p$ and $q$ i.e. that vanishes on $p$ but not on $q$.
To see this observe that every section of $H^0(L)$ vanish on $p$, since $p=\bl$, but
there is a section $\sigma\in H^0(L)$ that does not vanish on $q$, otherwise $q\in\bl$.
Now consider $\sigma^{2g+1}\in H^0((2g+1)L)$, this section separates $p$ from $q$.
\end{proof}

Theorem \ref{t5} is as best as possible as proved in the following:
\begin{propo}
Given a polarized manifold $(M,L)$ with 
$\Delta(M,L) \ = \ d(M,L) \ = \ 1$, sectionally hyperelliptic of
type $(-)$, the least integer $t$ such that $tL$ is very ample is $2g+1$ i.e.
$K_M+rL$ is not very ample for $r\leq n+1$.
\end{propo}
\begin{proof}
By contradiction suppose that $2gL$ is very ample and consider a regular ladder of $(M,L)$. 
Since the ladder is regular, the very ampleness of $2gL$ implies the very ampleness of the line bundles $L_{M_i}$.\\
In particular we have that $2g p$ is very ample, where $p=\bl$
and $2g p$ is a divisor on $C$, where $C$ is the one-dimensional element of the ladder.
Note that $C$ is an hyperelliptic curve and that $g(M,L)=g(C)$.\\
Now we may compute canonical bundle by adjunction.
By \ref{f1} we have $K_M=(2g-n-1)L$, then we obtain:
$K_D=(K_M+L)_D=(2g-n)L_D$, and by induction: $K_C=(2g-2)L_C$. 
But $L_C=\bl=p$ so we have $K_C=(g-1)2p$.
Since $C$ is hyperelliptic we know that $h^0(2p)=2$.
Now consider the canonical map ${\varphi}_{\mid K_C\mid}$ di $C$ and let
$<\sigma,\tau>$ is a basis of $H^0(2p)$, a basis of $H^0(K_C)$ is:
$<{\sigma}^{g-1},{\sigma}^{g-2}{\tau},\ldots,{\tau}^{g-1}>$. \\
Now observe that $<\sigma^g,\sigma^{g-1}\tau,\ldots,\tau^g>$ is a basis of $H^0(K_C+2p)$. 
In fact by Riemann-Roch we have $h^0(K_C+2p)=g+1$. 
This theorem implies that also $\varphi_{\mid K_C+2p\mid}$ is a double covering
of ${\p}^1$  followed  by the Veronese embedding,
but this implies that $(2g+1)L_C=K_C+2p$ is not very ample.
\end{proof}
\subsection{Surfaces of general type}
Let $S$ be a minimal surface of general type (i.e. of Kodaira dimension 2) and
let $K_S$ be its canonical bundle.
It is known \cite{bo} that $K$ is ample outside the set of $(-2)$-curves, i.e. $KC=0$ if and only if $C$ is a $(-2)$-curve.\\
In literature there are results \cite{bo} on the $C$-isomorphism property 
of the maps $\varphi_n=\varphi_{\mid nL\mid}$,
there are also results \cite{ci} on the generation of the canonical ring:      

\[
\bigoplus_{t\geq 1}H^0(S,tK_S).
\]

As theorem \ref{tm} shows there must be a connection between these results.
In the following table the results on the generation of the canonical ring, 
for surfaces with $q=H^1(S,{\mathcal O}_S)=0$, are compared with the very 
ampleness of $rK_S$. For each pair of values of $K_S^2$ and 
$P_g(S)=h^0(S,K_S)$ are compared the $C$-isomorphism property of 
$\varphi_{\mid rK_S\mid}$ and the $r$-generation of $K_S$.
\begin{center}
\begin{tabular}{||c|c|c|c||} \hline
$K_S^2$   & $\varphi{rK_S}$ is a $C$-isomorphism & $P_g$ & $K_S$ is
$r$-generated\\ \hline
        &                      &   $0$ & $\leq 6$\\ \cline{3-4} 
  $1$   &           $5$        &   $1$ & $\leq 4$ \\ \cline{3-4}
        &                      &   $2$ & $5$\\ \hline                
        &                      &   $0$ & $\leq 6$ \\ \cline{3-4}
        &                      &   $1$ & $\leq 4$ \\ \cline{3-4}           
  $2$   &           $4$        &   $2$ & $3 (\star)$\\ \cline{3-4}
        &                      &   $3$ & $4$\\ \hline    
        & $3K_s$ is spanned    &   $0$ & $\leq 6$\\ \cline{3-4}
        & and $\varphi_3$ is a &   $1$ & $\leq 4$\\ \cline{3-4}
  $3$   & birational morphism  &   $2$ & $3 (\star)$\\ \cline{3-4}      
        &                      &   $3$ & $2 (\bullet)$\\ \hline
\end{tabular}
\end{center}
Where $(\star)$ means that the result is true if $\mid K_S\mid$
contains an irreducible element and $(\bullet)$ means that the result is true if 
$K_S$ is base point free and $\varphi_{\mid K_S\mid}(S)$ is not a rational surface of degree
$r-1$ in $\p^r$ ($r=h^0(K_S)$).\\
We observe that $r$-generation of $K_S$ is not so far from the
very ampleness of $rK_S$ and by theorem $\ref{tm}$ we have also the 
following propositions:

\begin{propo}
\label{p1}
Let $S$ be a surface of general type with $P_g=K^2=2$ and $q=0$. 
If the generic element of $\mid K_S\mid$ is irreducible, then 
$\varphi_{\mid 3K_S\mid}$ is a $C$-isomorphism outside $\bs\mid K_S\mid$.
\end{propo}
\begin{proof}
By theorem $3.7$ of \cite{ci}, $K_S$ is $3$-generated so the proof follows 
from theorem \ref{tm}.
\end{proof}

similarly one can prove:
\begin{propo}
Let $S$ a surface of general type with $q=0$, $P_g\geq 3$,  
$K_S^2\geq 3$ and $K_S$ spanned, then we have two cases:
\begin{itemize}
\item  $\varphi_{\mid K_S\mid}(S)$ is a rational surface of degree
$r-1$ in $\p^r$ ($r=h^0(K_S)$)
\item $2K_S$ is a $C-isomorphism$
\end{itemize}
\end{propo}
\subsection{Del Pezzo surfaces}
In the following table is resumed the correspondence between very 
ampleness and $r$-beneration foa a Del Pezzo Surface. Here $L=-K_S$ 
\begin{center} \begin{tabular}{||c|c|c||}\hline
degree & very ampleness &  $r$-generation\\ \hline

$L^2=3$   & $L$  & $1$-generated\\ \hline
$L^2=2$   & $2L$ & $2$-generated\\ \hline
$L^2=1$   & $3L$ & $3$-generated\\ \hline
\end{tabular}
\end{center}

\bibliographystyle{plain}

\end{document}